# Modified Meshless Methods based on the Moving least Squares Method for Solution the Model of the Risk of Microcephaly Induced by the Zika Virus (ZIKV)

M. Matin far$^{a*}$, M. Pourabd$^a$, E. Taghizade$^a$

a. Department of Mathematics, Science of Mathematics Faculty, University of Mazandaran, Iran

**Abstract**

The aim of this work is the application of the Meshfree methods for solving systems of stiff ordinary differential equation. This method is based on the Moving least squares (MLS), generalized moving least squares (GMLS) approximation and Modified Moving least squares (MMLS) method. GMLS makes a considerable reduction in the cost of numerical methods. In fact, GMLS method is effect operator on the basis polynomial rather than the complicated MLS shape functions. Besides that the modified MMLS approximation method avoids undue a singular moment matrix. This allows the base functions to be of order greater than two with the same size of the support domain as the linear base functions. We also show the estimation of the error propagation obtained of the numerical solution of the systems of stiff ordinary differential equation. Some examples are provided to show that the GMLS and MMLS methods are more reliable (accurate) than classic MLS method.Finally, the (our) proposed methods are validated by solving ZIKV model which is a system of ODEs.

**keywords**: Moving least squares, Modified Moving least squares,General Moving least squares, systems of stiff ordinary equations , Model of the risk of microcephaly induced by the Zika virus (ZIKV), Error estimate

## 1 Introduction

The numerical solution of large-scale scientific problems expressed as the stiff systems of ordinary differential equations is now a significant issue in the areas of chemical engineering, nonlinear mechanics, biochemistry and life sciences. In computational mathematics, improving and extending methods for the solution of these systems of functional equations is one of the attractive areas for research [7–12]. Additionally, Haar wavelets method [1, 2], Adomian decomposition method [3], the variational iteration method [4], analytical Backlund transformation method [5] and semi-analytic methods , HPM [6] are other approaches proposed. The Moving Least Squares (MLS) method is an efficient numerical meshless method. The flexibility and high accuracy of the MLS approximation method are the main points of its development by different authors and use of it for solving a vast number of problems. Lancaster and Salkauskas [17] have introduced the Moving Least Squares (MLS) approximation method inspired by researchers McLain [14, 15] and Shepard [16] . Also, the MLS approximation has significant applications in different problems of the computational mathematics that were cited in [13, 17–20, 22, 24–30]. In this paper, is proposed a new modification of meshless numerical approximation method based on moving least squares(MLS) approximation that was introduced by G.R.Joldesa and et. [21]. This modifies allows, quadratic base functions ($m = 2$) to be used with the same size of the support domain as linear base functions ($m = 1$). In fact, the main property of this modification is that prevent the singular moment matrix in the MLS based methods. The paper presents the results of the solution of the stiff system of ordinary differential equations to the modified MLS method.

---

*Corresponding author. Address:Department of Mathematics, University of Mazandran, P.O.Box 47415-95447, Babolsar,Iran. Fax:01135302460. E-mail address:m.matinfar@umz.ac.ir

In [38], the relation of generalized moving least squares (GMLS) approximation and Backus-Gilbert optimality is investigated and then an application to numerical integration was performed in [39]. Moreover, in [40] the error estimates for generalized moving least squares (GMLS) derivatives approximations of a Sobolev function in $L_p$ norms and extends them for local weak forms of DMLPG methods is presented. It should be noted that, GMLS derivatives approximations are much easier to evaluate at the considerably lower cost in comparison with the MLS derivatives approximations.

The rest of this paper is organized as follows. In section 2 the concepts of three methods classic MLS , GMLS, and MMLS were reviewed and in section 3 stiff system of ordinary differential equations is presented and was implemented the proposed methods on it. In section 4, we investigated error bounds of the linear system. In section 5 by some numerical experiments, the error estimates and convergence rates were vindicated.. finally, in the last section, the proposed methods was applied for solving the ZIKV problem which is a system of ODE.

## 2 Outline of the method

### 2.1 MLS approximation method and GMLS

consider a function $U \in \mathcal{F}(U)$ where $\mathcal{F}(U)$ is a function space with certain smoothness (usually a Sobolev space) in abounded domain $\Omega$. To implementation, Assume that the nodal points be as a set of N points such as the following ,

$$X = \{\mathbf{x}_1, \mathbf{x}_2, \ldots, \mathbf{x}_N\}$$

where $X \subseteq \Omega \subseteq \mathbb{R}^d$ and $\Omega$ be a nonempty and bounded set.

Let $\mathbb{P}_m^d$ the space of d-variable polynomial of degree at most $m, m \in \mathbb{N}$, and let $B_m = \{p_1, p_2, \ldots p_m\}$ be any basis of $\mathbb{P}_m$. Also, $\Omega_{\mathbf{x}} = B_r(\mathbf{x}) \cap \Omega$ where the function is approximated at $\mathbf{x} \in X$ and it is the neighborhood of any points and the domain of definition of the MLS approximation, which are located in the problem domain $\Omega$.

To approximate the unknown vector of functions $U(\mathbf{x})$ for any of $u_i, i = 1, 2, \ldots, n$ in $\Omega_{\mathbf{x}}$, a number of randomly located nodes are selected, $x_i, i = 1, 2 \ldots, N$, and then the MLS approximate $u_i^h(\mathbf{x})$ of $u_i(\mathbf{x}) \in U(\mathbf{x}), \forall \mathbf{x} \in \Omega_{\mathbf{x}}$, can be defined as:

$$u_i^h(\mathbf{x}) = \sum_{j=1}^{m} a_j(\mathbf{x}) p_j(\mathbf{x}), i = 1, 2, \ldots n \tag{1}$$

where $a_j, j = 1, 2, \ldots, m$ are MLS shape functions which are chosen such that:

$$J(\mathbf{x}) = [P.\mathbf{a} - \mathbf{u}_i]^T.W.[P.\mathbf{a} - \mathbf{u}_i], i = 1, 2, \ldots, n \tag{2}$$

is minimized. where $W(x)$ is the diagonal matrix carrying the weights $w_i(\mathbf{x})$ on its diagonal, with $w_i(x) > 0$, the matrices $P$ are defined as, $P = [\mathbf{p}^T(\mathbf{x}_1), \mathbf{p}^T(\mathbf{x}_2), \ldots \mathbf{p}^T(\mathbf{x}_N)]_{N \times (m+1)}^T$ and $\mathbf{u}^h = [u_1^h, u_2^h, \ldots u_n^h]$. It is important to note that $u_i, i = 1, 2, \ldots n$, are the artificial nodal values, and not the nodal values of the unknown trial function $u^h(\mathbf{x})$ in general. With respect to $\mathbf{a}(\mathbf{x})$ and $u_i$ will be obtained,

$$\mathbf{a}(\mathbf{x}) = W(\mathbf{x}) P^T (P^T W(\mathbf{x}) P)^{-1} \mathbf{p}(\mathbf{x}), \tag{3}$$

the matrices $A(\mathbf{x}) = W(\mathbf{x}) P^T$ is symmetric positive definite matrix for all $\mathbf{x} \in \Omega$ and $\mathbf{p} = (p_1, p_2, \ldots, p_m)^T$. With computing $\mathbf{a}(\mathbf{x})$, $u_i^h$ can be obtained as follows:

$$u_i^h(\mathbf{x}) = \sum_{j=1}^{N} a_j(\mathbf{x}) u_i(\mathbf{x}_j) = \mathbf{a}^T . \mathbf{u}_i \tag{4}$$

As we know, each redial basis function that define in [23] can be used as weight function, we can define $w_j(r) = \phi(\frac{r}{\delta})$ where $r = ||\mathbf{x} - \mathbf{x}_i||_2$ (the Euclidean distance between $\mathbf{x}$ and $\mathbf{x}_j$) and $\phi : \mathbb{R}^d \longrightarrow \mathbb{R}$ is a nonnegative function with compact support. In this paper, we will use following weight functions and will compare them to each other, corresponding to the node $j$, in the numerical examples.

a: Guass weight function

$$w(r) = \begin{cases} \frac{\exp(\frac{-r^2}{c^2}) - \exp(\frac{-\delta^2}{c^2})}{1 - \exp(\frac{-\delta^2}{c^2})} & 0 \leq r \leq \delta \\ 0 & elsewhere. \end{cases} \tag{5}$$

b: RBF weight function

$$w = \begin{cases} (1-r)^6 (6 + 36r + 82r^2 + 72r^3 + 30r^4 + 5r^5) & 0 \leq r \leq \delta \\ 0 & elsewhere. \end{cases} \tag{6}$$

c: Spline weight function

$$w = \begin{cases} 1 - 6(\frac{r}{\delta})^2 + 8(\frac{r}{\delta})^3 - 3(\frac{r}{\delta})^4 & 0 \leq r \leq \delta \\ 0 & elsewhere. \end{cases} \tag{7}$$

Where $c$ is constant and is called shape parameter,Also $\delta$ is the size of support domain. If, further, $\phi$ is sufficiently smooth, derivatives of $U$ are usually approximated by derivatives of $U^h$,

$$D^\alpha u_i \mathbf{x} \approx D^\alpha u_i^h(\mathbf{x}) = \sum_{j=1}^{N} D^\alpha a_j(\mathbf{x}) u_i(\mathbf{x}_j), \mathbf{x} \in \Omega \tag{8}$$

More generally, consider a function $U \in \mathcal{F}(U)$ and fix a functional $\lambda \in \mathcal{F}(U)^*$, the dual space of $\mathcal{F}(U)$. In MLS,for each $_i \in U)$ , $\lambda(u_i)$ can be approximated by [36]

$$\widehat{\lambda(u_i)} \approx \lambda(\widehat{u_i}) = \sum_{j=1}^{N} \lambda(a_j) u_i(x_j)$$

where the functional $\lambda$ can be for instance point evaluations, derivative or integral operators, and etc. so it requires the act of functional $\lambda$ on shape functions $a_j$ and sometimes needs many calculations, especially when $\lambda$ is a complicated functional. The GMLS approximation can be specified as follows

$$\lambda(\hat{u}_i) = \sum_{j=1}^{N} (a_{\lambda,j}) u_i(x_j)$$

where $\lambda(u_i)$ directly approximated from nodal values $u_i(\mathbf{x}_1), \ldots, u_i(\mathbf{x}_N)$, and $a_{\lambda,j}$ are functions correspond to functional $\lambda$. The notion of GMLS almost similar the classic MLS method, so for the finite space $\mathbb{P}_m = span\{p_1, p_2, \ldots, p_m\}$, i.e.

$$\sum_{j=1}^{N} a_{j,\lambda} p(\mathbf{x}_j) = \lambda(p), \quad p \in P.$$

Then GMLS approximation $\widehat{\lambda(u_i)}$ to $\lambda(u_i)$ can be obtained as $\widehat{\lambda(u_i)} = \lambda(p^*)$, where $p^* \in P$ is minimizing the weighted least squares error functional

$$\sum_{j=1}^{N} w_j [p(\mathbf{x}_j) - u_i(\mathbf{x}_j)]^2, \tag{9}$$

among all $p \in P$ [37], so optimal $a_\lambda^T$ can also be determined as follows,

$$a_\lambda^T = \lambda(p^T) A^{-1}(\mathbf{x}) B(\mathbf{x}), \tag{10}$$

where $A$ and $B$ are defined as before and

$$\lambda(p^T) = [\lambda(p_1), \ldots, \lambda(p_m)] \in R^m$$

In fact, the main idea of GMLS approximation is effect operator on the basis polynomial functions.

## 2.2 Modified MLS approximation method

One of the common problems in Classic MLS method is the singularity of the moment matrix A in irregularity chosen nodal points. To avoid the nodal configurations which lead to a singular moment matrix, the usual solution is to enlarge the support domains of any nodal point. But it isn't an appropriate solution, in [21] to tackle such problems is proposed a modified Moving least squares(MMLS)approximation method. This modifies allows, base functions $m \geq 2$ to be used with the same size of the support domain as linear base functions ($m = 1$). We should note that,impose additional terms based on the coefficients of the polynomial base functions is the main view of the modified technique . As we know, in the basis function $\mathbf{p}(\mathbf{x})$ is

$$\mathbf{p}(\mathbf{x}) = [1, x, x^2, \ldots, x^m]^T \tag{11}$$

where $\mathbf{x} \in \mathbb{R}$, Also the correspond coefficients $a_j$, that should be determined are [28]:

$$\mathbf{a}(\mathbf{x}) = [a_1, a_x, a_{x^2}, \ldots, a_{x^m}]^T \tag{12}$$

For obtaining these coefficients, the functional (2) rewrite as follows:

$$\overline{J}(\mathbf{x}) = \sum_{j=1}^{m}(\mathbf{P}^T(\mathbf{x}_j)\mathbf{a}(\mathbf{x}) - u_i(\mathbf{x}_j))^2 w_i(\mathbf{x}) + \sum_{\nu=1}^{m-2} \overline{w}_\nu(\mathbf{x})\overline{\mathbf{a}}_\nu^2(\mathbf{x}), i = 1, 2, \ldots, n \tag{13}$$

Where $\overline{w}$ is a vector of positive weights for the additional constraints, also $\overline{\mathbf{a}} = [a_{x^2}, a_{x^3}, \ldots, a_{x^m}]^T$ is the modified matrix.
The matrix form of (2) is as follows:

$$\overline{J}(\mathbf{x}) = [P.\mathbf{a} - \mathbf{u}_i]^T.W.[P.\mathbf{a} - \mathbf{u}_i] + \mathbf{a}^T H \mathbf{a}, i = 1, 2, \ldots, n \tag{14}$$

where $H$ is as,

$$H = \begin{bmatrix} O_{2,2} & O_{m-2,m-2} \\ O_{2,2} & diag(\overline{w}) \end{bmatrix}, \tag{15}$$

where, $O_{i,j}$ is the null matrix. By minimizing the functional (14), the coefficients $a(\mathbf{x})$ will be obtained. So we have

$$\overline{A}(\mathbf{x})\mathbf{a}(\mathbf{x}) = B(\mathbf{x})\mathbf{u}_i, \tag{16}$$

where

$$\overline{A} = P^T.W.P + H \tag{17}$$

And the matrics $B(\mathbf{x})$ is determined as the same before. So we have

$$\varphi_m(\mathbf{x}) = \mathbf{a}(\mathbf{x}) = \mathbf{p}^T(\mathbf{x})\overline{A}^{-1}(\mathbf{x})B(\mathbf{x}) \tag{18}$$

where $\varphi_m(\mathbf{x})$ is the shape function of the MMLS approximation method.

## 3  Stiff Systems of Ordinary Differential Equations

In this section, we use MLS approximation method for numerical solution of the Stiff system of ordinary differential equations so consider the following differential equation

$$A(U) - F(\mathbf{x}) = 0, \ U(0) = U_0, \ \mathbf{x} \in \Omega \tag{19}$$

with boundary conditions,

$$B(U, \frac{\partial U}{\partial \mathbf{x}}) = 0, \ \mathbf{x} \in \partial\Omega.$$

where A is a general differential operator, $U_0$ is an initial approximation of Eq. (19), $F(\mathbf{x})$ is a vector of known analytical functions on the domain $\Omega$ and $\partial\Omega$ is the boundary of $\Omega$. The operator can be divided by $A = L + N$, where $L$ is the linear part, and $N$ is the nonlinear part of its. So Eq. (19) can be, rewritten as follows;

$$L(U) + N(U) - F(\mathbf{x}) = 0 \tag{20}$$

We assume that $\mathbf{a} = \{a_1, a_2, \ldots, a_m\}$ are the MLS shape functions so in order to solve system (20), N nodal points $x_i$ are selected on the $\Omega$, which $\{x_i | i = 1, 2, \ldots, N\}$ is q-unisolvent. The distribution of nodes could be selected regularly or randomly. Then instead of $u_j$ from $U$, we can replace $u_j^h$ from Eq.(4). So we have

$$u_j^h(\mathbf{x}) = \sum_{i=1}^{N} a_i(\mathbf{x}) u_j(\mathbf{x}_i) \tag{21}$$

where $j = 1, 2, \ldots, n$ is the number of unknown functions. we estimate the unknown functions $u_i$ by Eq.(21), so the system (20) becomes the following form

$$L(\sum_{i=1}^{N} a_i(\mathbf{x})u_1(\mathbf{x}_i), \sum_{i=1}^{N} a_i(\mathbf{x})u_2(\mathbf{x}_i), \ldots, \sum_{i=1}^{N} a_i(\mathbf{x})u_n(\mathbf{x}_i)) + \tag{22}$$
$$N(\sum_{i=1}^{N} a_i(\mathbf{x})u_1(\mathbf{x}_i), \sum_{i=1}^{N} a_i(\mathbf{x})u_2(\mathbf{x}_i), \ldots, \sum_{i=1}^{N} a_i(\mathbf{x})u_n(\mathbf{x}_i)) = (f_1(\mathbf{x}), f_2(\mathbf{x}), \ldots, f_n(\mathbf{x})) + \mathbf{r}(\mathbf{x}).$$

where $\mathbf{r}(\mathbf{x})$ is residual error function which vanishes to zero in collocation points thus by installing the collocation points $\mathbf{y}_r; r = 1, 2, \ldots, N$, so

$$L(\sum_{i=1}^{N} a_i(\mathbf{y}_r)u_1(\mathbf{x}_i), \sum_{i=1}^{N} a_i(\mathbf{y}_r)u_2(\mathbf{x}_i), \ldots, \sum_{i=1}^{N} a_i(\mathbf{y}_r)u_n(\mathbf{x}_i)) + \tag{23}$$
$$N(\sum_{i=1}^{N} a_i(\mathbf{y}_r)u_1(\mathbf{x}_i), \sum_{i=1}^{N} a_i(\mathbf{y}_r)u_2(\mathbf{x}_i), \ldots, \sum_{i=1}^{N} a_i(\mathbf{y}_r)u_n(\mathbf{x}_i)) =$$
$$\sum_{i=1}^{N} L(a_i(\mathbf{y}_r))u_1(\mathbf{x}_i), \sum_{i=1}^{N} L(a_i(\mathbf{y}_r))u_2(\mathbf{x}_i), \ldots, \sum_{i=1}^{N} L(a_i(\mathbf{y}_r))u_n(\mathbf{x}_i)) +$$

$$N(\sum_{i=1}^{N} a_i(\mathbf{y}_r)u_1(\mathbf{x}_i), \sum_{i=1}^{N} a_i(\mathbf{y}_r)u_2(\mathbf{x}_i), \ldots, \sum_{i=1}^{N} a_i(\mathbf{y}_r)u_n(\mathbf{x}_i)) =$$
$$(f_1(\mathbf{y}_r), f_2(\mathbf{y}_r), \ldots, f_n(\mathbf{y}_r))$$

therefore

$$CU = \begin{bmatrix} L(a_1(y_1)) & L(a_2(y_1)) & \ldots & L(a_N(y_1)) \\ L(a_1(y_2)) & L(a_2(y_2)) & \ldots & L(a_N(y_2)) \\ \vdots & & & \\ L(a_1(y_N)) & L(a_2(y_N)) & \ldots & L(a_N(y_N)) \end{bmatrix} \begin{bmatrix} u_1(x_1) & u_2(x_1) & \ldots & u_n(x_1) \\ u_1(x_2) & u_2(x_2) & \ldots & u_n(x_2) \\ \vdots & & & \\ u_1(x_N) & u_2(x_N) & \ldots & u_n(x_N) \end{bmatrix} \quad (24)$$

And the matrix form of (23) as follows

$$C_{N\times N}U_{N\times n} + N_{N\times n}(\mathbf{a}, U) = F_{N\times n}(\mathbf{y}_r) \quad (25)$$

where

$$\begin{aligned} C_i &= [L(a_1(\mathbf{y}_r)), \ldots, L(a_N(\mathbf{y}_r))]_{i=1}^{n} \\ U_i &= [(u_i(x_1), u_i(x_2), \ldots, u_i(x_N))^T]_{i=1}^{n} \\ F(\mathbf{y}_r) &= ([(f_1(\mathbf{y}_r))_{r=1}^N]^T, [(f_2(\mathbf{y}_r))_{r=1}^N]^T, \ldots [(f_n(\mathbf{y}_r))_{r=1}^N]^T)^T. \end{aligned} \quad (26)$$

by imposing the initial conditions at $t = 0$, we have

$$(\sum_{i=1}^{N} a_i(0)u_1(t_i), \sum_{i=1}^{N} a_i(0)u_2(t_i), \ldots, \sum_{i=1}^{N} a_i(0)u_n(t_i)) = U_0 \quad (27)$$

and Solving the non-linear system (25) and (27), lead to finding quantities $u_j(x_i)$. Then the values of $u_j(x)$ at any point $x \in \Omega$, can be approximated by Eq.(4) as:

$$u_j(\mathbf{x}) \simeq \sum_{i=1}^{N} a_i(\mathbf{x})u_j(\mathbf{x}_i)$$

## 4 Bound Error

This section devoted to the error estimation. In [32], has obtained error estimates for moving least square approximations in the one-dimensional case. Also in [33], is developed for functional in n-dimensional and was used the error estimates to prove an error estimate in Galerkin coercive problems. In this work, have improved error estimate for the systems of stiff ordinary differential equations.

Given $\delta > 0$ let $W_\delta \geq 0$ be a function such that $supp(w_\delta) \subset \overline{B_\delta(0)} = \{z||z| \leq \delta\}$ and $X_\delta = \{x_1, x_2, \ldots, x_n\}$, $n = n(\delta)$, a set of points in $\Omega \subset \mathbb{R}$ an open interval and let $U = (u_1, u_2, \ldots, u_N)$ be the unknown function such that $u_{i1}, u_{i2}, \ldots, u_{in}$ be the values of the function $u_i$ in those points, i.e., $u_{i,j} = u_i(x_j), i = 1, \ldots, N, j = 1, \ldots, n$. A class of functions $W = \{\omega_j\}_{j=1}^{N}$ is called a partition of unity subordinated to the open covering $I_N$ if it possesses

the following properties:

- $W_j \in C_s^0, s > 0 \text{ or } s = \infty,$
- $\sup(\omega_j) \subseteq \bar{\Lambda}_j,$
- $\omega_j(x) > 0, x \in \Lambda_j,$
- $\sum_{i=1}^{N} \omega_j = 1 \text{ for every } x \in \bar{\Omega}$

There is no unique way to build a partition of unity as defined above [34]
As we know, the MLS approximation is well defined if we have a unique solution at every point $x \in \bar{\Omega}$. for minimization problem. And non-singularity of matrix $A(x)$, ensuring it is .In [33] the error estimate was obtained with the following assumption about the system of nodes and weight functions $\{\Theta_N, W_N\}$ :

We define
$$\langle u, v \rangle = \sum_{j=1}^{n} w(x - x_j) u(x_j) v(x_j)$$

then
$$\|u\|_x^2 = \sum_{j=1}^{n} w(x - x_j) u(x_j)^2$$

Also for vector of unknown functions, we define

$$\|U\|_\infty = max\{|u_i|_x, i = 1, 2, \ldots, N\}$$

are the discrete norm on the polynomial space $\mathbb{P}_m^1$ if the weight function $w$ satisfy the following properties.

**a**: For each $x \in \Omega$, $w(x - x_j) > 0$ at least for $(m+1)$ indices $j$.

**b**: For any $x \in \Omega$, the moment matrix $A(x) = w(x) P^T$ is nonsingular.

**Definition 4.1.** *Given* $\boldsymbol{x} \in \bar{\Omega}$, *the set* $ST(\boldsymbol{x}) = \{j : \omega_j \neq 0\}$ *will be called the star of* $x$.

**Theorem 4.1.** *( [?, ?].) A necessary condition for the satisfaction of Property **b** is that for any* $\boldsymbol{x} \in \bar{\Omega}$

$$n = card(ST(\boldsymbol{x})) \geq card(\mathbb{P}_m) = m + 1$$

For a sample point $\mathbf{c} \in \bar{\Omega}$, if $ST(\mathbf{c}) = \{j_1, \ldots j_k\}$, the mesh-size of the star $ST(\mathbf{c})$ defined by the number is $h(ST(\mathbf{c})) = \max\{h_{j1}, \ldots h_{jk}\}$.

**Assumptions.** Consider the following global assumptions about parameters. There exist

($a_1$) An over bound of the overlap of clouds:
$$E = \sup_{c \in \bar{\Omega}}\{card(ST(c))\}.$$

($a_2$) Upper bounds of the condition number:
$$CB_q = \sup_{c \in \bar{\Omega}}\{CN_q(ST(c)), q = 1, 2\}.$$

where the numbers $CN_q(ST(\mathbf{c}))$ are computable measures of the quality of the star $ST(c)$ which defined in Theorem7 of [18].

($a_3$) An upper bound of the mesh-size of stars:
$$R = \sup_{c \in \bar{\Omega}}(hST(c)).$$

($a_4$) An uniform bound of the derivatives of $\{w_j\}$. That is the constant $G_q > 0, q = 1, 2$, such that
$$\|D^\mu W_j\|_{L_\infty} \leq \frac{G_q}{R^{|\mu|}} \quad 1 < \mu < q,$$

($a_5$) There exist the number $\gamma > 0$ such that any two points $\mathbf{x}, \mathbf{y} \in \Omega$ can be joined by a rectifiable curve $\Gamma$ in $\Omega$ with length $|\Gamma| \leq \gamma \|\mathbf{x}\text{-}\mathbf{y}\|$.
Assuming all these conditions, Zuppa [33] proved.

**Lemma 4.1.** $U = (u_1, u_2, \ldots u_n)$ such that $u_i \in C^{m+1}(\bar{\Omega})$ and $\|U\|_\infty = u_k$, $1 < k < n$, There exist constants $C_q, q = 1 or 2$,

$$C_1 = C_1(\gamma, d, E, G_1, CB_1),$$
$$C_2 = C_1(\gamma, d, E, G_2, CB_1, CB_2),$$

then
$$\left\|D^\mu U - D^\mu U^h\right\|_\infty < C_q R^{q+1-|\mu|} \|u_k^{(m+1)}\|_{L^\infty(\Omega)} \quad 0 < \mu < q \tag{28}$$

**Proof:** see [15]

## 4.1 System of ODE

If in Eq.(20) the non-linear operator $N$ be zero, we have

$$L(U) = (f_1, f_2, \ldots, f_n) \tag{29}$$

where $U$ is the vector of unknown function and $L$ is a matrix of derivative operators,

$$L(U(.)) = \sum_{\varsigma=1}^{n} \lambda_\varsigma \frac{\partial^\varsigma}{(.)^\varsigma} U(.). \tag{30}$$

And from (Eq.21), we define

$$U^h(t) = \sum_{i=1}^{N} a_i(t) U(t_i)$$

where $(a_i)_{i=1}^{N}$ are the MLs shape functions defined on the interval $[0,1]$ satisfying the homogeneous counterparts of the boundary conditions in Eq. (19). Also if the weight function $w$ possesses $k$ continuous derivatives then the shape functions $a_j$ is also in $C^k$ [20]. By the collocation method, is obtained an approximate solution $U^h(t)$. And demand that

$$L^h(U(.)) = \sum_{\varsigma=0}^{n} \lambda_\varsigma \frac{\partial^\varsigma}{(.)^\varsigma} U^h(.) \qquad (31)$$

where $(\lambda = 0 \, or \, 1)$. It is assumed that in the system of ODE derivative of order at most $n = 2$. Each of the basis functions $a_i$ has compact support contained in $(0,1)$ then the matrix $C$ in Eq.(26) is a bounded matrix. If $\delta$ be fixed, independent of $N$, then the resulting system of linear equations can be solved in $O(N)$ arithmetic operations.

**Lemma 4.2.** *Let $U = (u_1, u_2, \ldots u_n)$ and $F = (f_1, f_2, \ldots f_n)$ so that $u_i \in C^{m+1}(\overline{\Omega})$ $m \geq 1$ and $\| u_i \|_\infty = u_k, k \in \{1, 2, \ldots, n\}$ where $\Omega$ be a closed, bounded set in $R$. Assume the quadrature scheme is convergent for all continuous functions on $\Omega$. Further, assume that the stiff system of ODE 19 with the fixd initial condition is uniquely solvable for given $f_i \in C(\Omega)$. Moreover take a suitable approximation $U^h$ of $U$ Then for all sufficiently large $n$, the approximate matrix $L$ for linearly case exist and are uniformly bounded, $|L| \leq M$ with a suitable constant $M < \infty$. For the equations $L(U) = F$ and $L^h(U) = F$ we have*

$$E_t = \|L(U(t)) - L^h(U(t))\|_\infty$$

*so that*

$$\|E_t\|_\infty \leq C_q K(\lambda, \varsigma) R^{m+1-\mu} \|u_k^{(m+1)}\|_{L_\infty}.$$

**Proof.** we have

$$\| L(U(t)) - L^h(U(t)) \|_\infty = \| \sum_{\varsigma=0}^{n} \lambda_\varsigma \frac{\partial^\varsigma}{t^\varsigma} U(t) - \sum_{\varsigma=0}^{n} \lambda_\varsigma \frac{\partial^\varsigma}{t^\varsigma} U^h(t) \|_\infty$$

so due to the lemma(4.1),

$$\begin{aligned}
\| L(U(t)) - L^h(U(t)) \|_\infty &\leq \sum_{\varsigma=0}^{n} |\lambda_\varsigma| \, \| \frac{\partial^\varsigma}{t^\varsigma} U(t) - \frac{\partial^\varsigma}{t^\varsigma} U^h(t) \|_\infty \\
&\leq \max_i \sum_{\varsigma=0}^{n} |\lambda_\varsigma| | \frac{\partial^\varsigma}{t^\varsigma} u_i(t) - \frac{\partial^\varsigma}{t^\varsigma} u_i^h(t) | \\
&\leq \sum_{\varsigma=0}^{n} C_q |\lambda_\varsigma| \|u_k^{(m+1)}\|_{L_\infty} R^{m+1-\varsigma}
\end{aligned}$$

where should be $m \geq \varsigma$ so,

$$\sum_{\varsigma=0}^{n} |\lambda_\varsigma| R^{m+1-\varsigma} \leq K(\lambda,\varsigma) R^{m+1-\mu}$$

where $\mu$ is the highest order derivative And $K(\lambda,\varsigma) = \sum_{\varsigma=0}^{n} |\lambda_\varsigma|$, so demanded that

$$\|E_t\|_\infty \leq C_q K(\lambda,\varsigma) R^{m+1-\mu} \|u_k^{(m+1)}\|_{L_\infty}.$$

It should be noted that in the nonlinear system the upper bound of error depends on the nonlinear operator.

## 5 Numerical Examples

In this section, we apply MLS and MMLS methods to solve the following test examples.

### 5.1 Example 1

Consider the following nonlinear stiff systems of ODEs [6]

$$\begin{cases} u_1'(t) = -1002 u_1(t) + 1000 u_2^2(t) \\ u_2'(t) = u_1(t) - u_2(t) - u_2^2(t) \end{cases}$$

With the initial condition $u_1(0) = 1$ and $u_2(0) = 1$. The exact solution is

$$u_1(t) = exp(-2t)$$

$$u_2(t) = exp(-t)$$

. In this numerical example, two scheme are compared and as explained the main task of the modified method tackle the singularity of the moment matrix. Table1 presents the maximum relative error by MLS on a set of evaluation points(with $h = 0.1 and 0.02$) and $\delta = 4h and 3h$. Also in Table2. MLS and MMLS at different number of nodes for $h = 0.004$ and $\delta = 5h and 8h$, were compared.

Table 1: Maximum relative errors by MLS , Example1.

|  | m=2, $\delta$ =4h | | m=2, $\delta$ =3h | |
| --- | --- | --- | --- | --- |
| h | $u_1$ | $u_2$ | $u_1$ | $u_2$ |
| 0.1 | $5 \times 10^{-3}$ | $4.1 \times 10^{-4}$ | $8.85 \times 10^{-4}$ | $2.2 \times 10^{-3}$ |
| 0.02 | $5.8 \times 10^{-2}$ | $6.5 \times 10^{-5}$ | $5.42 \times 10^{-4}$ | $6.52 \times 10^{-5}$ |

Table 2: Maximum relative errors for h=0.004 by MMLS and MLS , Example1.

|      | m=3,$\delta$ =5h | | m=3, $\delta$ =8h | |
|------|------|------|------|------|
| Type | $u_1$ | $u_2$ | $u_1$ | $u_2$ |
| MLS  | $1.03 \times 10^0$ | $0.98 \times 10^1$ | $1.01 \times 10^0$ | $9.2 \times 10^0$ |
| MMLS | $9.23 \times 10^{-4}$ | $9.22 \times 10^{-4}$ | $6.89 \times 10^{-4}$ | $6.96 \times 10^{-4}$ |

## 5.2 Example 2

in this example, we consider $U(t) = (\frac{1}{47}(95\exp^{(-2t)} - 48exp(-96t)), \frac{1}{47}(48\exp(-96t) - \exp(-2t)))$ as the exact solution and $U(0,0) = (1,1)$ as the initial conditions for the following system of ODE,

$$\begin{cases} x'(t) = -x(t) + 95y(t) \\ y'(t) = -x(t) - 97y(t) \end{cases}$$

Table3 presents the maximum relative norm of the errors on a fine set of evaluation points (with $h = 0.004$) and $\delta = 5h$ for MLS and GMLS at different type of weight functions. As seen in this table, one major advantage of GMLS is that the computational time used by GMLS is less than MLS.

Table 3: Maximum relative errors by MLS $t \in [0,5]$ ,$h = 0.004$, , Example2.

|  | MLS | | | GMLS | | |
|---|---|---|---|---|---|---|
| weight type | $u_1$ | $u_2$ | Cputime | $u_1$ | $u_2$ | Cputime |
| Guass | $3.06 \times 10^{-3}$ | $9.92 \times 10^{-5}$ | 61.1706 | $8.5 \times 10^{-4}$ | $6.5 \times 10^{-4}$ | 0.5598 |
| Spline | $5.06 \times 10^{-4}$ | $5.3 \times 10^{-4}$ | 64.5897 | $1.93 \times 10^{-2}$ | $4.23 \times 10^{-4}$ | 0.6714 |
| RBF | $6.407 \times 10^{-4}$ | $3.02 \times 10^{-4}$ | 59.1790 | $6.9 \times 10^{-3}$ | $6.9 \times 10^{-3}$ | 0.5768 |

# 6 Application

In this section, is investigated a mathematical model based on the system of ODE equations, which trust into a dynamical system of nonlinear differential equations to account for the risk of microcephaly incidence caused by the Zika virus. in [35], is introduced a theoretical model to describe the dynamics of the population incidence of the infected pregnant women that might display fetal microcephaly induced by the ZIKV virus. The variables of the model are as follows:

-$x_1(t)$: average number of susceptible people,

-$x_2(t)$ : average number of ZIKV infected pregnant women that may induce fetal microcephaly,

-$x_3(t)$: average number of persons infected by ZIKV,
-$y_1(t)$: average number of non-carrier mosquitoes,
-$y_2(t)$: average number of virus-carrier mosquitoes,
we set $X(t)$ as the total people population and $Y(t)$ total population of mosquitoes at time $t$.
Also according to research cited, the parameters used for the simulations are,
$\eta$: constant flux of susceptible people,
$\mu$: the people natural death rate, $\beta$: the virus transmission probability from the virus carrier mosquitoes to the susceptible people,
$\sigma$: the virus transmission probability from the infected pregnant women to the non-carrier mosquitoes,
$\gamma$: virus transmission probability from infected people to the non-carrier mosquitoes,
$\epsilon$: the adult mosquitoes death rate,
$\theta$: the recovery rate of the infected pregnant women, $\alpha$: the infected people recovery rate,
$f$ : the fraction of infected people,
$1 - f$: is the fraction of pregnant women infected by ZIKV.
So by the above variables and parameters, the differential equations system of the infectious process will be shown as following,

$$\begin{aligned}
\frac{\partial x_1(t)}{\partial t} &= \eta - \beta \frac{y_2(t)}{Y(t)} x_1(t) - \mu x_1(t) \\
\frac{\partial x_2(t)}{\partial t} &= (1-f)\beta \frac{y_2(t)}{Y(t)} x_1(t) - (\theta + \mu) x_2(t) \\
\frac{\partial x_3(t)}{\partial t} &= f\beta \frac{y_2(t)}{Y(t)} x_1(t) - (\alpha + \mu) x_3(t) \\
\frac{\partial y_1(t)}{\partial t} &= \rho - \sigma \frac{x_2(t)}{X(t)} y_1(t) - \gamma \frac{x_3(t)}{X(t)} y_1(t) - \epsilon y_1(t) \\
\frac{\partial y_2(t)}{\partial t} &= \sigma \frac{x_2(t)}{X(t)} y_1(t) + \gamma \frac{x_3(t)}{X(t)} y_1(t) - \epsilon y_2(t)
\end{aligned} \quad (32)$$

where $\eta, \alpha, \mu, \theta, \rho, \epsilon > 0 \, and \, \beta, \sigma, \gamma, f \in (0,1)$ and initial conditions are $x_1(0) = x_{10}, x_2(0) = x_{20}, x_3(0) = x_{30}, y_1(0) = y_{10}, y_2(0) = y_{20}$. Taking into account the importance of these previous reports about of the relation between ZIKV and microcephaly in newborns, solving this system, is more crucial than ever. We solving the system was accomplished by using the data of Table 4. which have been previously reported in [35]. Computing the unknown functions of the system of ZIKV by MLS and MMLS is the main goal. These systems of ODE correspond to the solutions of each population, so by analyzing the model, we have considered the initial conditions as $X(0) = (0,0,0)$ and $Y(0) = (0,0)$. Also, in computations we put $h = 5, \delta = 5h$ and $t_i \in [0,40]$. In Figure1. compared the behavior of the susceptible peoples , infected preg-nant women and persons infected by ZIKV fo different value of $f = 0.3 and 0.6$. Also, over about 20 days the populations tend to stabilize. Figure 2, investigate the effect of the fraction of the infected people $f$ on infection process of pregnant women population.

Table 4: Parameters Values

| Parameter | $\gamma$ | $\beta$ | $\sigma$ | $\epsilon$ | $\alpha$ | $\mu$ | $\theta$ | $\rho$ | $\eta$ | $f$ |
|---|---|---|---|---|---|---|---|---|---|---|
| Value | 0.773 | 0.7913 | 0.6 | 0.0352 | 0.14 | 0.0003 | 0.05 | 30 | 20 | 0.3, 0.6 |

Figure 1: Behavior of the $x_1$,$x_2$ and $x_3$, the left is $f = 0.3$ and right is $f = 0.6$ by MLS

Figure 2: Behavior of the infected pregnant women$x_2$ by the zika virus
$\langle - \bullet -f = 0.3, -o - f = 0.6, - + -f = 0.85, - * -f = 1 \rangle$

Figure 3: Behaviorof the infected pregnant women$x_2$ respect to the number of virus-carrier mosquitoes $y_2$

Finally, Figure 3, Indicative the number of virus-carrier mosquitoes effect on the pregnant women that risk acquiring induce fetal microcephaly.

Figure 4: Behavior of the $x_1$,$x_2$ and $x_3$, the left is $f = 0.3$ and right is $f = 0.6$ by MMLS

In Modified moving least squares scheme, we put $\mu = 0.1$ and all the other value of the parameter such that $\delta$, $m$ and $h$ for computing the shape function are the same as MLS. Also, the spline weight function is used. The results are shown in Figure 4. It should be noted that by using the proposed method possibility of the investigating the system of the infection process by ZIKV of each population has been prepared. And according to the results of the implementation of the method on the system of ODE, these results is trusted.

# 7 Conclusion

In this paper, three meshless techniques called moving least squares, modified Moving least squares, and generalized moving least squares approximation are applied for solving the system of stiff ordinary differential equations. comparing the results obtained by these methods with the results obtained by the exact solution shows that the moving least squares methods are the reliable and accurate methods for solving stiff problems. Meshless methods are free from choosing the domain and this makes it suitable to study real- world problems. We solve the problem that models the operation of Zika virus without any restriction on the domain (population under study) and the number of parameters and this is one of our achievements in this research.